\documentclass[preprint,12pt,authoryear]{elsarticle}
\usepackage[english]{babel}
\usepackage{enumitem}
\usepackage{graphicx}
\usepackage{dcolumn}
\usepackage{bm}
\usepackage{amsmath,amsfonts}
\usepackage{amssymb}
\usepackage{amsthm}
\usepackage{physics}
\usepackage{color}
\usepackage{xcolor}
\usepackage[final,colorlinks,linkcolor=blue,anchorcolor=blue,urlcolor=blue,citecolor=blue]{hyperref}
\usepackage{caption}
\usepackage{subcaption}
\usepackage{booktabs}
\usepackage{siunitx}
\usepackage{epstopdf}
\usepackage{diagbox}
\usepackage{framed}
\usepackage{dsfont}
\usepackage{float}
\usepackage[margin=2.5cm]{geometry}
\usepackage[normalem]{ulem}

\newtheorem{proposition}{Proposition}

\begin{document}

\title{Elementary Tail Bounds on the Hypergeometric Distribution}

\author[1,2,3]{Vaisakh Mannalath$^{\ast}$}
\ead{vmannalath@vqcc.uvigo.es}

\author[1,2,3]{Víctor Zapatero}

\author[1,2,3]{Marcos Curty}

\address[1]{Vigo Quantum Communication Center, University of Vigo, Vigo E-36310, Spain}
\address[2]{Escuela de Ingeniería de Telecomunicación, Department of Signal Theory and Communications, University of Vigo, Vigo E-36310, Spain}
\address[3]{AtlanTTic Research Center, University of Vigo, Vigo E-36310, Spain}

\begin{abstract}
We use  simple methods to derive three concentration bounds on the hypergeometric distribution. Comparison with existing results illustrates the advantage of these bounds across different regimes.
\end{abstract}
\maketitle
\begingroup
\renewcommand\thefootnote{}
\footnotetext{$^{\ast}$Corresponding author.}
\addtocounter{footnote}{0}
\endgroup
\section{Introduction}
As established by  \cite{Hoeffding1963} and  \cite{Chvátal}, the inequality due to ~\cite{Chernoff} for sums of independent random variables is also applicable to a hypergeometric random variable. Nevertheless, when applied to a standard hypergeometric distribution, $\text{Hypergeometric}(N, K, n)$, the Chernoff bound is insensitive to the sampling fraction, $n/N$, and consequently becomes increasingly looser as the latter increases. To mitigate this effect, one can resort to different exponential bounds carrying a dependence on $n/N$. As surveyed in~\cite{Greene-thesis}, these include a renowned tail inequality by ~\cite{Serfling} and subsequent results by ~\cite{Hush}, ~\cite{Chatterjee}, ~\cite{Rohde}, ~\cite{Goldstein}, ~\cite{Bardenet} and ~\cite{Greene}.

Here, we present an elementary alternative: by exploiting a basic permutation symmetry between the parameters $n$ and $K$ in the hypergeometric distribution, we ``reformulate" the Chernoff bound to retrieve the dependence on the sampling fraction. Despite the simplicity of this approach, it outperforms the standard Chernoff bound for $n>K$, and exhibits a clear advantage over the listed results in that regime.

More generally, the validity of the Chernoff inequality in the hypergeometric setting follows from the Poisson binomial representation of the hypergeometric distribution (\cite{Pitman,Hui}). This stronger result, combined with a well-known theorem by ~\cite{Hoeffding1956}, allows to obtain a tighter non-exponential bound based on the binomial cumulative mass function. Again, the bound is detached from $n/N$ and worsens as the latter increases, but for $n>K$ one can mitigate this effect by leveraging the permutation symmetry between $n$ and $K$.

Finally, we complement the binomial-based results by deriving a factorial moment bound for the hypergeometric distribution. In spite of the fact that these bounds are known to improve upon standard Chernoff estimates for non-negative integer-valued variables as shown in (\cite{Philips,Naveau}), to the best of our knowledge, the hypergeometric case that we address has not been considered in the literature before.

\section{Results}

Let $X \sim \text{Hypergeometric}(N, K, n)$ and $d \geq nK/N +1$. Then,
\begin{equation}
\label{eq:chernoff}
\Pr[X \geq d] \leq \exp\left[-n \operatorname{D}(d/n \|\, K/N)\right]
\end{equation}
by virtue of the Chernoff inequality, where
\begin{equation}
\operatorname{D}(x \,\|\, y) = x \ln\left(\frac{x}{y}\right) + (1 - x)\ln\left(\frac{1 - x}{1 - y}\right)
\end{equation}
is the Kullback-Leibler divergence between Bernoulli random variables of parameters $x$ and $y$.

On the other hand, it is known that $X$ can be represented as a sum of $n$ independent but non-identical Bernoulli variables (\cite{Hui}), sometimes referred to as a Poisson binomial. In other words, there exist $p_1,p_{2}\ldots p_n \in [0,1]$ such that
\begin{equation}
    X = \sum_{i=1}^{n} Y_{i}
\end{equation}
for $Y_{i}\sim \text{Bernoulli}(p_i)$. Note that, in particular, this implies that $\sum_{i}p_{i}=nK/N$. By virtue of this representation, a theorem by  \cite{Hoeffding1956} yields
\begin{equation}
\label{eq:beta}
\Pr[X \geq d] \leq I_{K/N}(d, n - d + 1)
\end{equation}
\\
for any integer $d\geq  n K/N +1$, where $I_x(a,b)$ denotes the regularized incomplete $\beta$ function (which characterizes the binomial cumulative mass function).

For a given ratio $K/N$ and  sample size $n$, both the Chernoff bound of Eq.~(\ref{eq:chernoff}) and the ``$\beta$ bound" of Eq.~(\ref{eq:beta}) are independent of the sampling fraction $n/N$, becoming loose as the fraction increases. This is easily understood in view of the departure of the hypergeometric distribution from the binomial distribution for large $n/N$, due to the impact that the no-replacement has on the statistics. Nevertheless, by leveraging a symmetry of the hypergeometric distribution, we obtain alternative bounds sensitive to $n/N$. To be precise, the variable $Z \sim \text{Hypergeometric}(N,n,K)$, which results from swapping $n$ and $K$ in the definition of $X$, fulfills $\Pr[Z=k]=\Pr[X=k]$ for all $k$~\footnote{That is to say, $X$ and $Z$ share the exact same probability mass function.}, such that 

\begin{equation}
\Pr[X\geq d]=\Pr[Z\geq d]
\end{equation}
for all $d$ as well. In particular, considering again an integer $d\geq  n K/N +1$, the Chernoff and $\beta$ bounds for $Z$ yield
\begin{equation}
\Pr[Z \geq d] \leq \exp\left[-K \operatorname{D}(d/K\|\, n/N)\right]
\end{equation}
and
\begin{equation}
\Pr[Z \geq d] \leq I_{n/N}(d, K - d + 1),
\end{equation}
respectively. Thus, we conclude the following.
\begin{proposition}
Let $X \sim \text{Hypergeometric}(N, K, n)$. Then,
\begin{equation}
\Pr[X \geq d] \leq \exp\left[-K \operatorname{D}(d/K \,||\, n/N)\right]
\end{equation}
and
\begin{equation}
\Pr[X \geq d] \leq I_{n/N}(d, K - d + 1)
\end{equation}
for all integer $d\geq{}nK/N+1$.
\end{proposition}
In short, the two binomial tail bounds for $Z$ are legitimate tail bounds for $X$ as well, with the added feature of being sensitive to the sampling fraction. Simulations indicate that exploiting the symmetry is advantageous in both bounds as long as $n>K$, for all integer $d\geq nK/N+1$. In other words, the preferred ``bounding binomial", between $\mathrm{Binomial}(n,K/N)$ and $\mathrm{Binomial}(K,n/N)$, seems to always be the one with the smallest number of trials~\footnote{This is to be expected, since the variance is smaller for the binomial with the smallest number of trials.}.

\begin{figure}[]
    \centering
    \begin{minipage}{0.49\textwidth}
        \includegraphics[width=\linewidth]{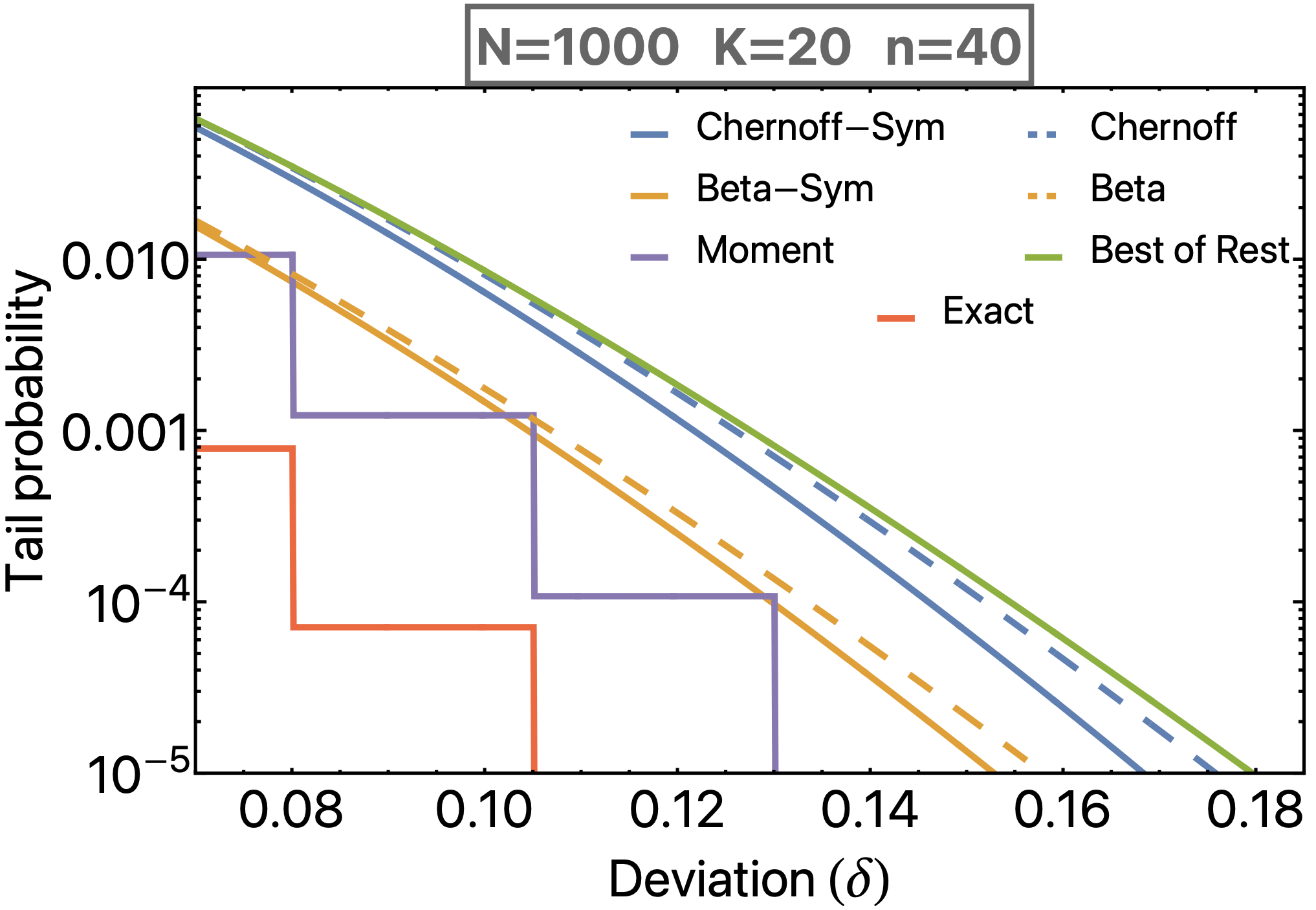}\vspace{0.3cm}
        \includegraphics[width=\linewidth]{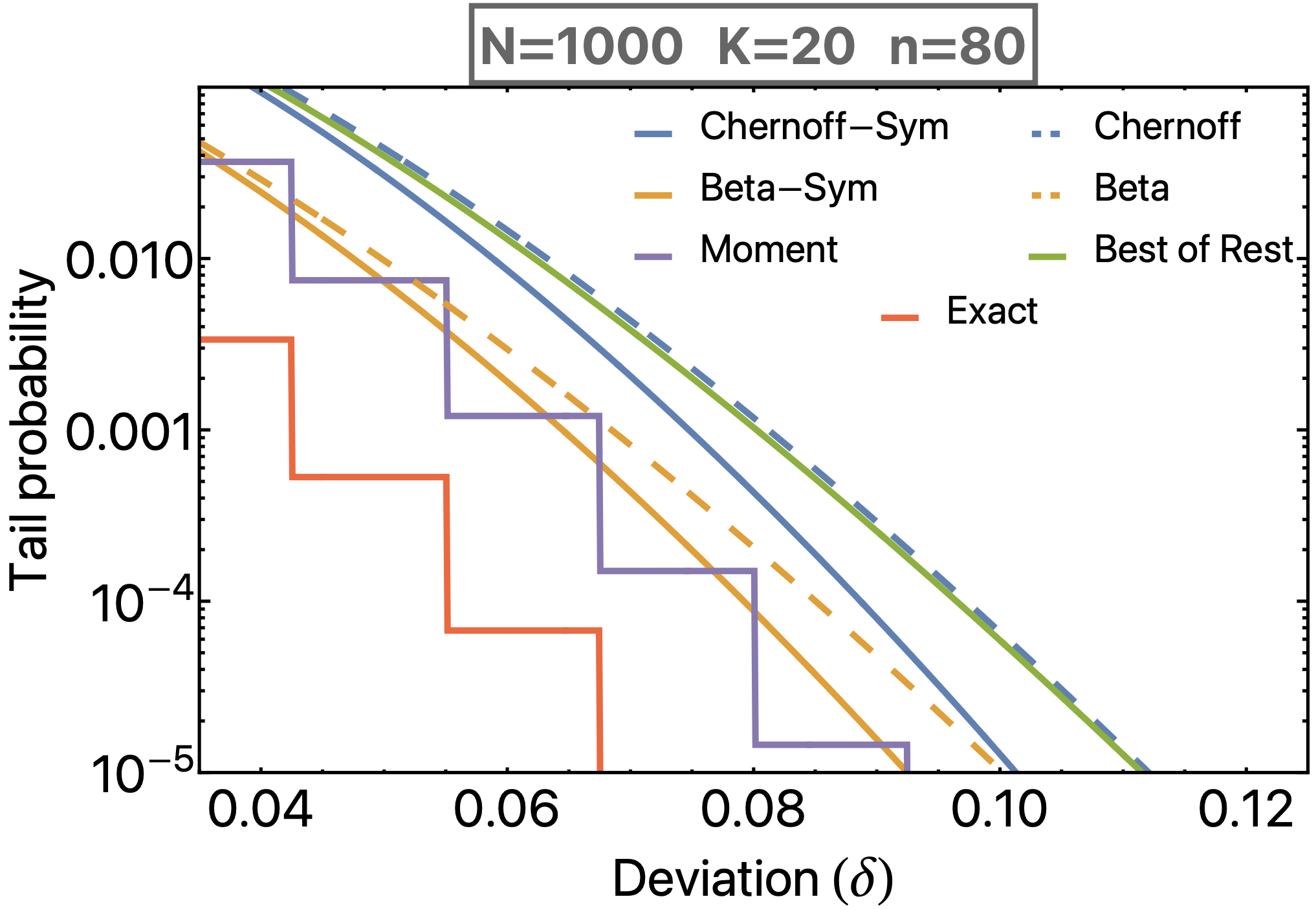}
    \end{minipage}
    \hfill
    \begin{minipage}{0.5\textwidth}
        \includegraphics[width=\linewidth]{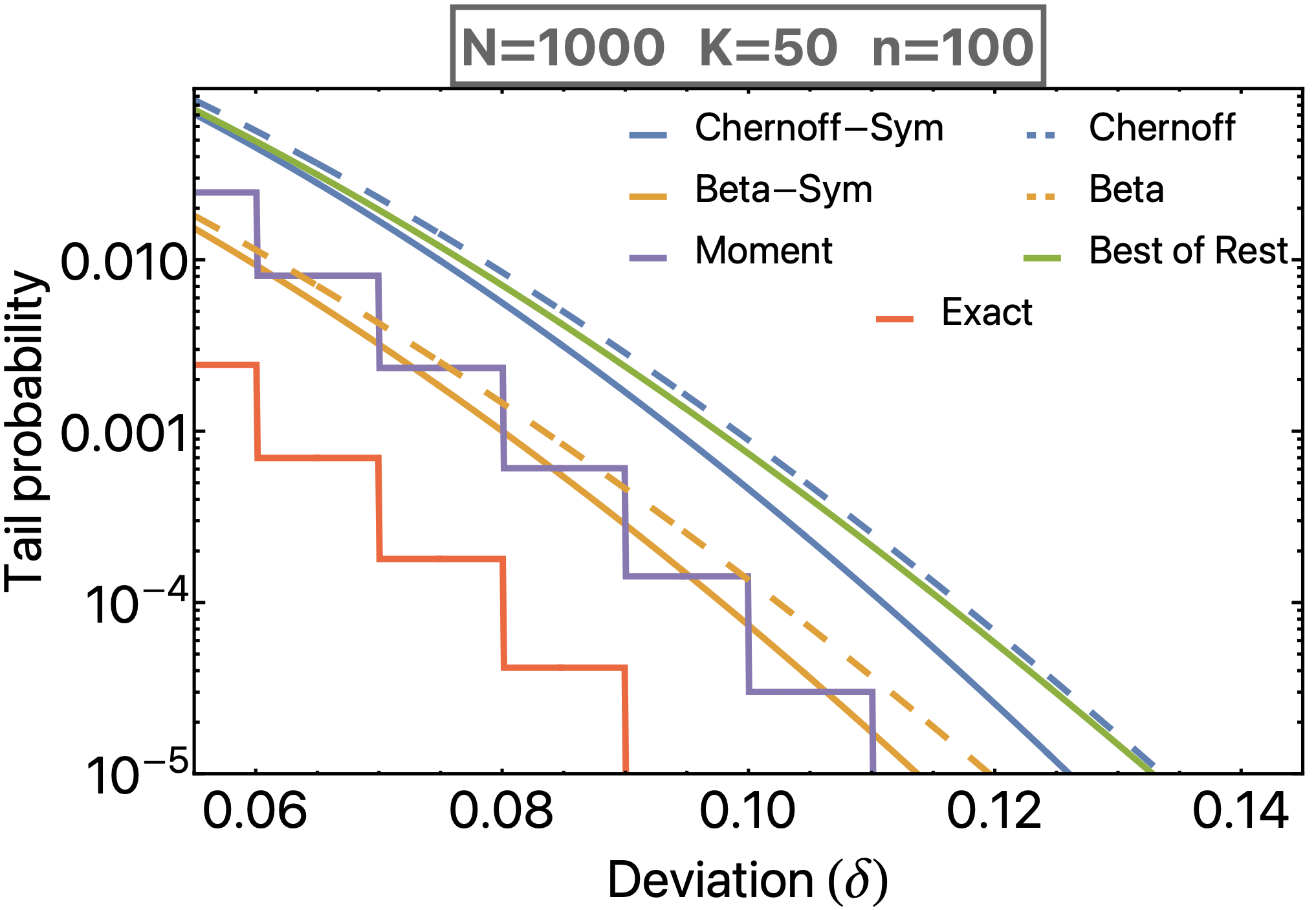}\vspace{0.3cm}
        \includegraphics[width=\linewidth]{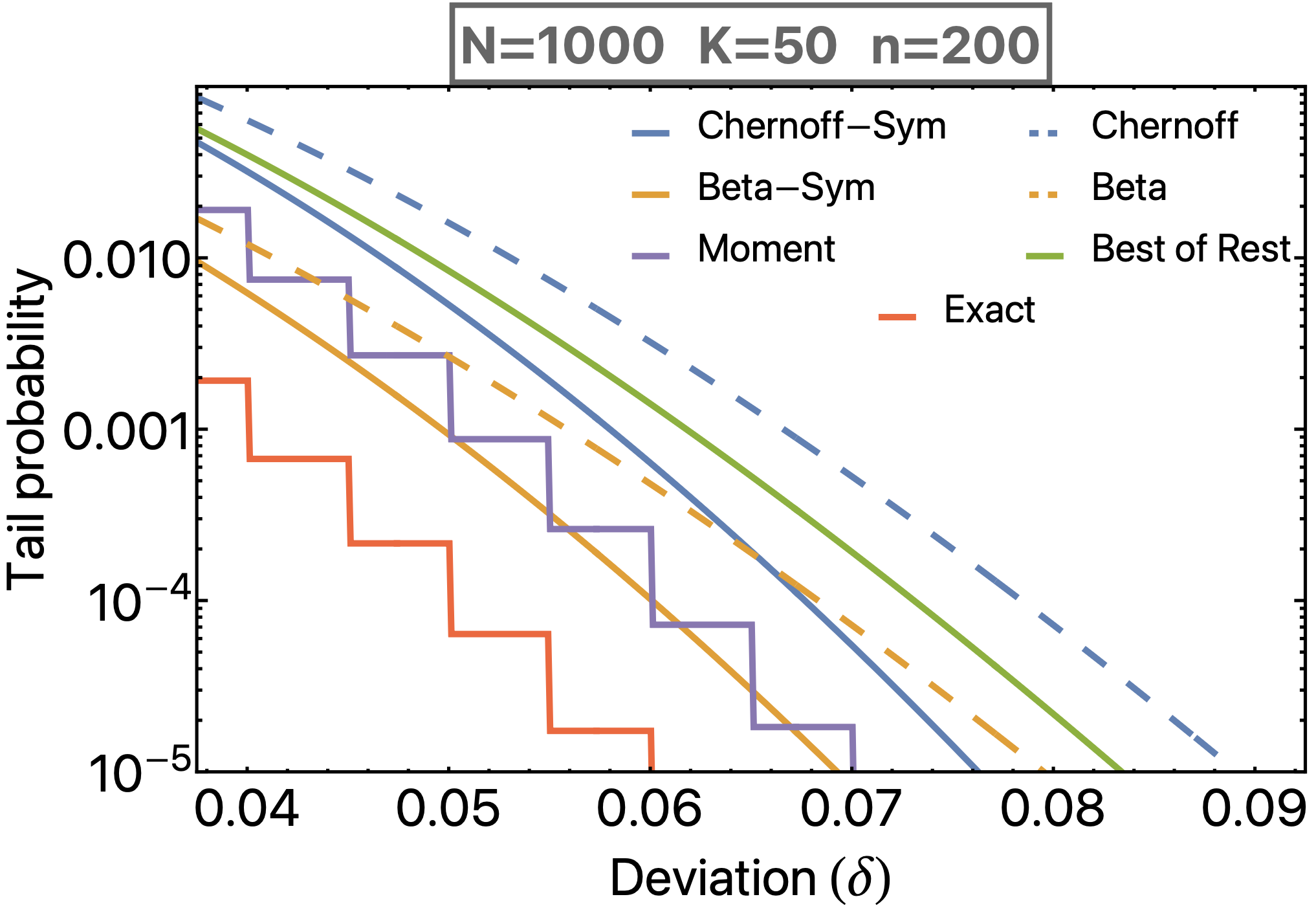}
    \end{minipage}
    \caption{Upper tail probability, $\Pr[X\geq n(K/N +\delta)]$, as a function of the relative deviation $\delta$ with respect to the mean, for $X\sim \text{Hypergeometric}(N,K,n)$ and a population size $N = 1000$. Left column: $K/N = 2\%$. Right column: $K/N = 5\%$. Blue and yellow lines: Chernoff and $\beta$ bounds with (solid) and without (dashed) invoking the permutation symmetry of $n$ and $K$. Green line: best-performing bound among ~\cite{Serfling}, ~\cite{Hush}, ~\cite{Chatterjee}, ~\cite{Goldstein}, ~\cite{Rohde},  ~\cite{Bardenet} and Theorem 1 and 3 of ~\cite{Greene}. Purple line: Factorial moment bound. Red line: exact cumulative mass function of $X$.}
    \label{fig:upper_tail_1000}
\end{figure}

\begin{figure}[]
    \centering
    \begin{minipage}{0.48\textwidth}
        \includegraphics[width=\linewidth]{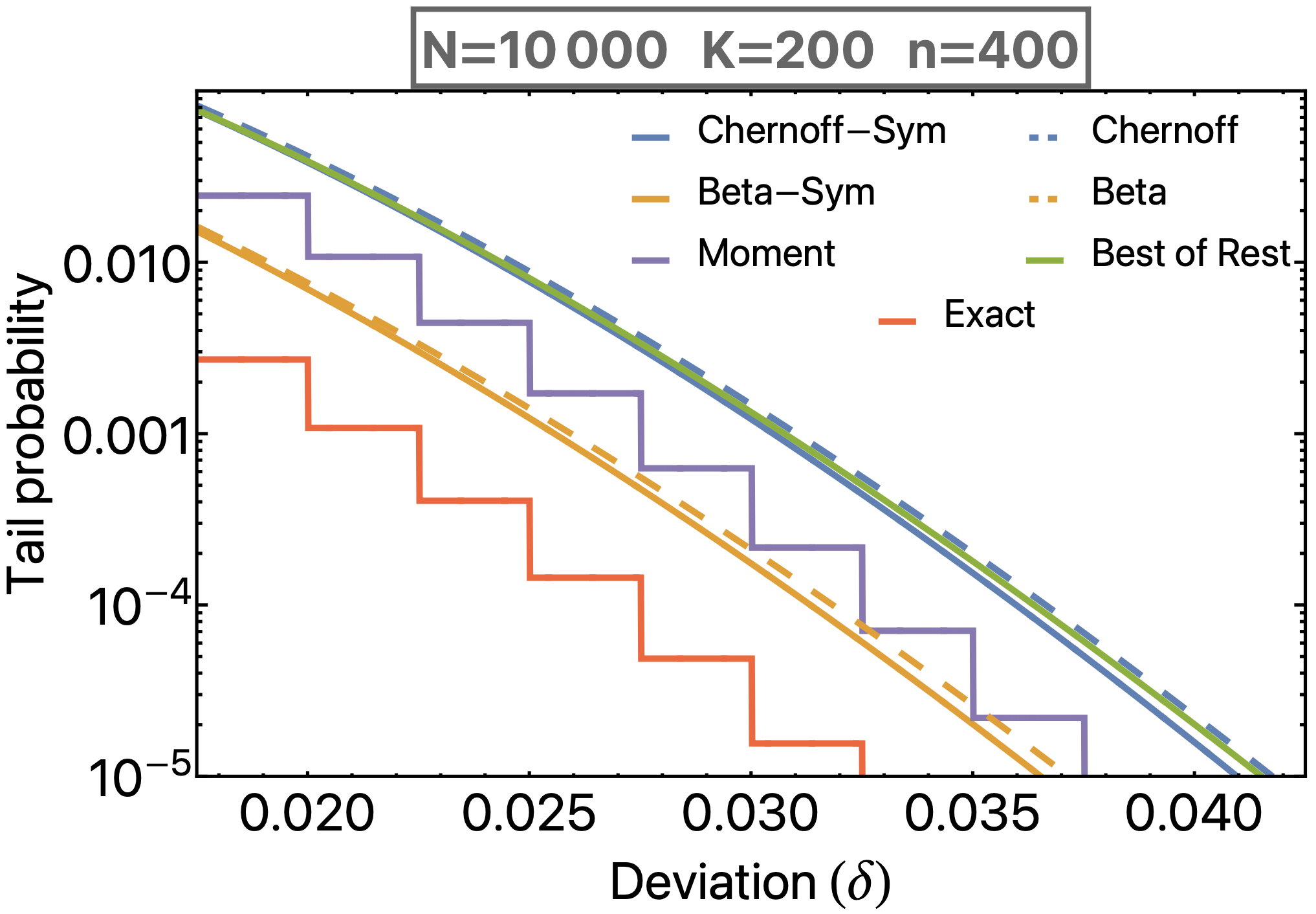}\vspace{0.3cm}
        \includegraphics[width=\linewidth]{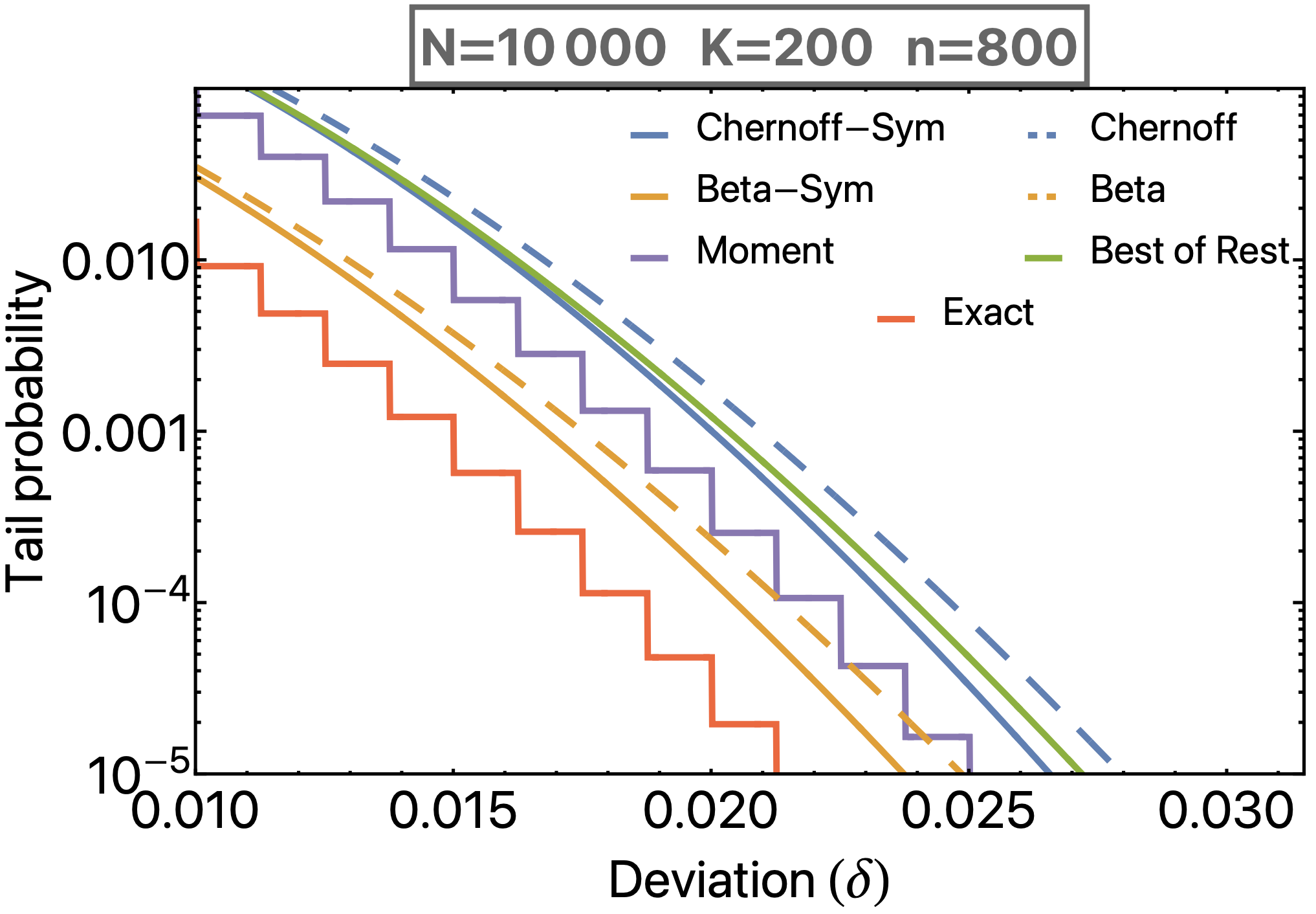}\vspace{0.3cm}
    \end{minipage}
    \hfill
    \begin{minipage}{0.48\textwidth}
        \includegraphics[width=\linewidth]{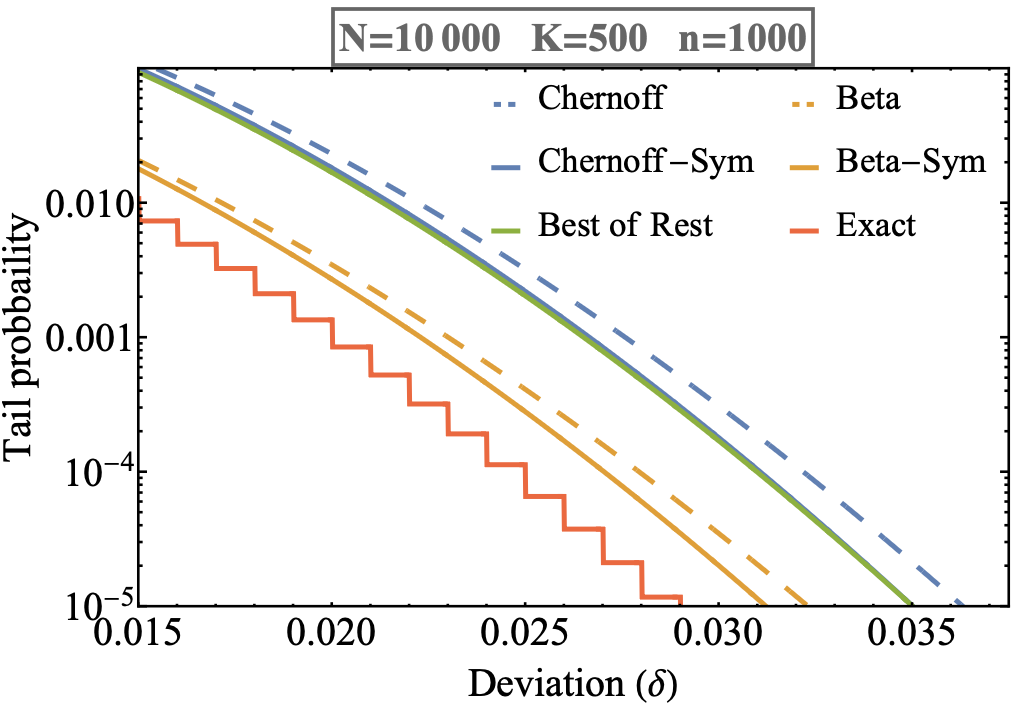}\vspace{0.3cm}
        \includegraphics[width=\linewidth]{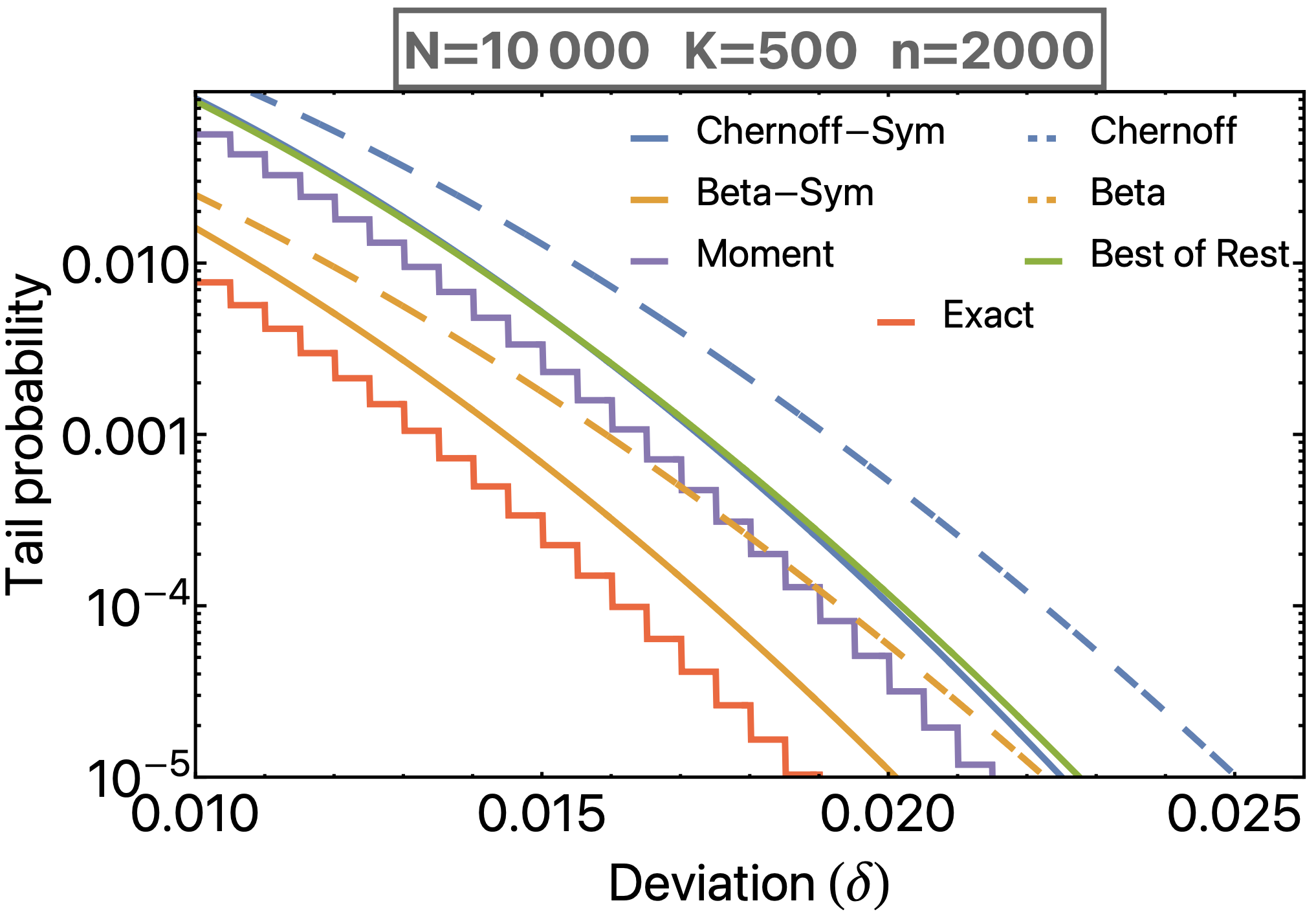}\vspace{0.3cm}
    \end{minipage}
    \caption{$\Pr[X\geq n(K/N +\delta)]$ as a function of $\delta$, for $X\sim \text{Hypergeometric}(N,K,n)$ and $N = 10000$. The followed criteria are common with those of Figure \ref{fig:upper_tail_1000}.}
    \label{fig:upper_tail_10000}
\end{figure}

In Figures~\ref{fig:upper_tail_1000} and \ref{fig:upper_tail_10000}, we illustrate this advantage with various examples, further including the exponential bounds of ~\cite{Serfling},~\cite{Hush}, ~\cite{Chatterjee}, ~\cite{Rohde}, ~\cite{Goldstein},  ~\cite{Bardenet} and  ~\cite{Greene}. For varying parameters, the plots display the upper tail probability $\Pr[X \geq n(K/N + \delta)]$ as a function of the relative deviation $\delta$. Figure \ref{fig:upper_tail_1000} corresponds to a fixed population size $N = 1000$, and Figure \ref{fig:upper_tail_10000} to $N = 10000$. Within each figure, the left column corresponds to $K/N = 2\%$, and the right column to $K/N = 5\%$.

Importantly, we note that additional symmetries of the hypergeometric distribution can be used to selectively improve the Chernoff and $\beta$ bounds in certain parameter regimes. Particularly, denoting $\Pr[X = k] = f(k; N, K, n)$, two other well-known symmetries to exploit are $f(k; N, K, n) = f(N-n-K+k; N, N-K, N - n)$ and $f(k; N, K, n) = f(N-n-K+k; N, N-n, N - K)$. These additional symmetries become useful when either $n$ or $K$ exceed $N/2$, simply because the Chernoff and $\beta$ bounds become tighter for the complement of the sample than for the sample itself (if $n>N/2$), and for the number of non-defective items rather than for the defective ones (if $K>N/2$).



Distinct from the symmetry-based approach discussed above, we derive a complementary result using factorial moment bounds. As discussed in the Introduction, the motivation is that this approach is guaranteed to provide strictly tighter estimates than the Chernoff bound ~(\cite{Philips,Naveau}).

By applying Markov's inequality to  a non-decreasing, non-negative function $\phi$  with $\phi(d)>0$,  we obtain

\begin{equation}
    \Pr[X\geq d]\leq \frac{\mathbb{E}[\phi(X)]}{\phi(d)}.
\end{equation}
Particularizing $\phi$ to the falling factorial with $t+1$ factors, $(x)_{t+1} = \prod_{i=0}^{t} (x - i)$, it follows that
\begin{equation}
    \Pr[X\geq d]\leq \inf_{0 \leq t<d}\frac{\mathbb{E}[X(X-1)\ldots(X-t)]}{d(d-1)\ldots (d-t)},
\end{equation}
and making use of the fact that (\cite{Potts})

\begin{equation}
    \mathbb{E}[(X)_{t+1}] = \frac{(K)_{t+1}(n)_{t+1}}{(N)_{t+1}},
\end{equation}
we obtain
\begin{equation}
\Pr[X \geq d] \leq \inf_{0 \leq t < d} \frac{(K)_{t+1}(n)_{t+1}}{(N)_{t+1}(d)_{t+1}}.
\end{equation}
Defining
\begin{equation}
f(t+1)=\frac{(K)_{t+1}(n)_{t+1}}{(N)_{t+1}(d)_{t+1}},
\end{equation}
we have that
\begin{equation}
\frac{f(t+1)}{f(t)}=\frac{(K - t)(n - t)}{(N - t)(d - t)},
\end{equation}
such that the condition $f(t+1)/f(t)<1$ translates into
\begin{equation}
t \leq \frac{N d - K n}{N + d - (K + n)}.
\end{equation}
This implies the following proposition.
\begin{proposition}
Let $X \sim \text{Hypergeometric}(N, K, n)$, and let $d\geq  n K/N +1$  be an integer. Then,
\begin{equation}
\Pr[X \geq d] \leq \frac{(K)_{t^* + 1}(n)_{t^* + 1}}{(N)_{t^* + 1}(d)_{t^* + 1}},
\end{equation}
where $ t^* = \left\lfloor \frac{N d - K n}{N + d - (K + n)} \right\rfloor$.
\end{proposition}

As anticipated, Figures~\ref{fig:upper_tail_1000} and \ref{fig:upper_tail_10000} show that the factorial moment bound is consistently tighter than the Chernoff bounds across the tested parameter regimes.

\newpage
\section*{Acknowledgements and Funding}
This work was supported by European Union’s Horizon Europe Framework Programme under the Marie Sklodowska-Curie Grant No. 101072637 (Project QSI), the Galician Regional Government (consolidation of research units: atlanTTic), the Spanish Ministry of Economy and Competitiveness (MINECO), the Fondo Europeo de Desarrollo Regional (FEDER) through the grant No. PID2024-162270OB-I00, MICIN with funding from the European Union NextGenerationEU (PRTRC17.I1) and the Galician Regional Government with own funding through the “Planes Complementarios de I+D+I con las Comunidades Autonomas” in Quantum Communication, the “Hub Nacional de Excelencia en Comunicaciones Cuanticas” funded by the Spanish Ministry for Digital Transformation and the Public Service and the European Union NextGenerationEU, the European Union’s Horizon Europe Framework Programme under the project Quantum Secure Networks Partnership (QSNP, grant agreement No 101114043) and the European Union via the European Health and Digital Executive Agency (HADEA) under the Project QuTechSpace (grant 101135225).

\section*{Availability of data and materials}
No datasets were generated or analysed during the current study.

\section*{Competing interests}
The authors declare no competing interests.

\section*{Authors' contributions}
V.M. and V.Z. developed the main ideas and performed the analysis of the proposed bounds. M.C. validated the results and refined the manuscript. All authors read and approved the final version.

%


\newpage
\begin{thebibliography}{20}
%
\bibitem[Hoeffding(1963)]{Hoeffding1963}
Hoeffding, W. Probability Inequalities for Sums of Bounded Random Variables. \textit{Journal of the American Statistical Association} \textbf{58}, 13–30 (1963).
%
\bibitem[Chvátal(1979)]{Chvátal}
Chvátal, V. The tail of the hypergeometric distribution. \textit{Discrete Mathematics} \textbf{25}, 285-287 (1979).
%
\bibitem[Chernoff(1952)]{Chernoff}
Chernoff, H. A measure of asymptotic efficiency for tests of a hypothesis based on the sum of observations. \textit{The Annals of Mathematical Statistics} \textbf{23}, 493-507 (1952).
%
\bibitem[Greene(2016)]{Greene-thesis}
Greene, E. Finite sampling exponential bounds (Doctoral dissertation, University of Washington, 2016).
%
\bibitem[Serfling(1974)]{Serfling}
Serfling, R. J. Probability inequalities for the sum in sampling without replacement. \textit{The Annals of Statistics} \textbf{2}, 39-48 (1974).
%
\bibitem[Hush and Scovel(2005)]{Hush}
Hush, D., \& Scovel, C. Concentration of the hypergeometric distribution. \textit{Statistics \& Probability Letters} \textbf{75}, 127-132 (2005).
%
\bibitem[Chatterjee(2007)]{Chatterjee}
Chatterjee, S. Stein’s method for concentration inequalities. \textit{Probability Theory and Related Fields} \textbf{138}, 305–321 (2007).
%
\bibitem[Rohde(2011)]{Rohde}
Rohde, A. Optimal calibration for multiple testing against local inhomogeneity in higher dimension. \textit{Probability theory and related fields} \textbf{149.3} 515-559 (2011).
%
\bibitem[Goldstein and Işlak(2014)]{Goldstein}
Goldstein, L., \& Işlak, Ü. Concentration inequalities via zero bias couplings. \textit{Statistics \& Probability Letters} \textbf{86}, 17-23 (2014).
%
\bibitem[Bardenet and Maillard(2015)]{Bardenet}
Bardenet, R., \& Maillard, O. A. Concentration inequalities for sampling without replacement. \textit{Bernoulli} \textbf{21,} 1361-1385 (2015).
%
\bibitem[Greene and Wellner(2017)]{Greene}
Greene, E., \& Wellner, J. A. Exponential bounds for the hypergeometric distribution. \textit{Bernoulli: Official Journal of the Bernoulli Society for Mathematical Statistics and Probability} \textbf{23}, 1911 (2017).
%
\bibitem[Pitman(1997)]{Pitman}
Pitman, J. Probabilistic bounds on the coefficients of polynomials with only real zeros. \textit{J. Comb. Theory Ser.
A} \textbf{77}, 279 (1997).
%
\bibitem[Hui and Park(2014)]{Hui}
Hui, S., \& Park, C. J. The representation of hypergeometric random variables using independent Bernoulli random variables. \textit{Communications in Statistics-Theory and Methods} \textbf{43}, 4103-4108 (2014).
%
\bibitem[Hoeffding(1956)]{Hoeffding1956}
Hoeffding, W. On the distribution of the number of successes in independent trials. \textit{The Annals of Mathematical Statistics}, 713-721 (1956).
%
\bibitem[Philips and Nelson(1995)]{Philips}
Philips, T. K., \& Nelson, R. The moment bound is tighter than Chernoff's bound for positive tail probabilities. \textit{The American Statistician} \textbf{49}, 175-178 (1995).
%
\bibitem[Naveau(1997)]{Naveau}
Naveau, P. Comparison between the Chernoff and factorial moment bounds for discrete random variables. \textit{The American Statistician} \textbf{51}, 40-41 (1997).
%
\bibitem[Potts(1953)]{Potts}
Potts, R. B. Note on the factorial moments of standard distributions. \textit{Australian Journal of Physics} \textbf{6}, 498-499 (1953).
%
\end{thebibliography}

\end{document}